\documentclass[a4paper,12pt]{article}
\usepackage{amssymb,amsmath,graphics,pstricks}

\newenvironment{dwd}{\par\noindent{\bf Proof.}}{\par\rightline{$\blacksquare$}}

\newenvironment{dwz}{\par\noindent{\bf Proof of Theorem 3.}}{\par\rightline{$\blacksquare$}}

\newtheorem{theo}{Theorem}
\newtheorem{lema}{Lemma}
\newtheorem{prop}{Proposition}
\newtheorem{defi}{Definition}

\newtheorem{rema}{Remark}
\newtheorem{colo}{Corollary}

\def\be#1\ee{\begin{equation}#1\end{equation}}
\newcommand{\ba}{\begin{eqnarray} }
\newcommand{\ea}{\end{eqnarray} }
\def\bt#1\et{\begin{theo}#1\end{theo}}
\def\bl#1\el{\begin{lema}#1\end{lema}}
\def\bp#1\ep{\begin{prop}#1\end{prop}}
\def\bd#1\ed{\begin{defi}#1\end{defi}}

\def\ccB{{\cal B}}

\def\ccM{{\cal M}}

\def\ccP{{\cal P}}

\def\va{\varepsilon}
\def\ra{\rightarrow}

\def\pint{-\hspace{-11pt}\int}

\def\E{\mathbf{E}}
\def\P{\mathbf{P}}

\def\N{{\mathbb N}}

\def\R{{\mathbb R}}

\def\ls{\leqslant}
\def\gs{\geqslant}
\def\for{\mbox{for}}
\def\For{\mbox{for all}}
\def\Bc{{B}}
\def\Bo{{B^{\circ}}}
\begin{document}
\title{\bf H{\"o}lder continuity of random processes}
\author{Witold Bednorz\footnotemark[1] \footnotemark[2] \footnotemark[3]}
\date{}
\maketitle
\footnotetext[1]{Department of Mathematic, University of Warsaw, Banacha 2, 02-097 Warsaw, Poland}
\footnotetext[2]{Partially supported by the Funds of Grant MENiN 1 P03A 01229}
\footnotetext[3]{E-mail: wbednorz@mimuw.edu.pl}

\begin{abstract}
For a Young function $\varphi$ and a Borel probability measure $m$ on a compact metric
space $(T,d)$ the minorizing metric is defined by
$$
\tau_{m,\varphi}(s,t):=\max\{\int^{d(s,t)}_0 \varphi^{-1}(\frac{1}{m(\Bc(s,\va))})d\va,
\int^{d(s,t)}_0 \varphi^{-1}(\frac{1}{m(\Bc(t,\va))})d\va\}.
$$
In the paper we extend the result of Kwapien and Rosinski \cite{Kwa}
relaxing the conditions on $\varphi$ under which there exists a constant $K$
such that
$$
\E\sup_{s,t\in T}\varphi(\frac{|X(s)-X(t)|}{K\tau_{m,\varphi}(s,t)})\ls 1,
$$
for each separable process $X(t)$, $t\in T$ which satisfies 
$\sup_{s,t\in T}\E \varphi(\frac{|X(s)-f(t)|}{d(s,t)})\ls 1$.
In the case of $\varphi_p(x)\equiv x^p$, $p\gs 1$
we obtain the somewhat weaker results. 

\end{abstract}

\section{Introduction}
Let $X$ be a topological space and $\ccB(X)$ its Borel $\sigma$-field. We denote by $\mathfrak{B}(X),\mathfrak{B}_b(X)$, $C(X),C_b(X)$
the set of all measurable, bounded measurable, continuous and bounded continuous functions respectively.
Furthermore $\ccP(X)$ denotes the family of all Borel, probability measures on $X$.
For each $\mu\in \ccP(X)$, $f\in\mathfrak{B}_b(X)$ and $A\in\ccB(X)$ we define
$$
\pint_A f(u)\mu (du):=\frac{1}{\mu(A)}\int_A f(u)\mu(du),
$$
where, we have used the convention $0/0=0$ (as we do throughout the whole paper). 
By $\mathrm{supp}(\mu)$ we denote the support of $\mu$.
\smallskip

\noindent
In the paper we consider finite Young functions; that is increasing convex functions 
$\varphi:[0,\infty)\ra [0,\infty)$ satisfying $\varphi(0)=0$, $\lim_{x\ra\infty}\varphi(x)=\infty$. For a simplicity we will be assuming also
that $\varphi(1)=1$. As in (\cite{Rao}, Def. 5, page 40), we let
$\triangle^2$ denote the set of all finite Young functions satisfying for some $c\gs 0, r>1$
{\renewcommand{\theequation}{$\triangle^2$}
\begin{equation}\label{pop}
\varphi(x)^2\ls \varphi(rx),\;\;\mbox{for some}\;\;\for\;\;x\gs c.
\end{equation}
and \addtocounter{equation}{-1}} let $\nabla'$ (see \cite{Rao}, Def 7, page 28) denote the set of all finite
Young functions $\varphi$ verifying for some $c\gs 0,r>1$ 
{\renewcommand{\theequation}{$\nabla'$}
\begin{equation}\label{pol}
\varphi(x)\varphi(y)\ls \varphi(rxy),\;\;\for\; x,y\gs c.
\end{equation}
Note \addtocounter{equation}{-1}} that if (\ref{pop}), resp. (\ref{pol}) holds for some $c>0$, then (\ref{pop}), resp. (\ref{pol}), holds for every $c'>0$ with appropriate choice of $r'$. If $h\in \mathfrak{B}(X)$ we let
$$
|h|^{\mu}_{\varphi}:=\inf\{a>0:\;\int_X\varphi(\frac{|h(s)|}{a})\mu(ds)\ls 1\},\;\;
\|h\|^{\mu}_{\varphi}:=\inf_{a>0}a(1+\int_X\varphi(\frac{|h(s)|}{a}))\mu(ds).
$$
denote the two Orlicz norms of $h$. Then $|\cdot|^{\mu}_{\varphi}$ and $\|\cdot\|^{\mu}_{\varphi}$ are semi-norms
on $\mathfrak{B}(X)$, satisfying $|h|^{\mu}_{\varphi}=0$ $\Leftrightarrow$ $\|h\|^{\mu}_{\varphi}=0$ $\Leftrightarrow$ $h=0$, $\mu$-a.e. 
Note that $|h|^{\mu}_{\infty}<\infty$ $\Leftrightarrow$ $\int_X \varphi(\frac{|h|}{a})<\infty$ for some $0<a<\infty$ 
$\Leftrightarrow$ $\|h\|^{\mu}_{\varphi}<\infty$ and recall that the Orlicz space $L^{\varphi}(\mu)$ is the set of all measurable functions
satisfying one of the three equivalent conditions (see \cite{Rao}). Then $(L^{\varphi}(\mu),|\cdot|_{\varphi})$ is
a complete semi-normed space. As we prove in Lemma \ref{1lem} semi-norms $|\cdot|^{\mu}_{\varphi}$ and $\|\cdot\|^{\mu}_{\varphi}$
are comparable.  
\smallskip

\noindent
Let $(T,d)$ be a fixed compact, metric space and
$m$ a fixed probability measure (defined on Borel subsets) on
$T$. For $x\in T$ and $\va\gs 0$, $\Bc(x,\va), \Bo(x,\va)$ denote respectively the closed and the open ball with the center at $x$ and the radius $\va$ i.e. $\Bc(x,\va)=\{y\in T:d(x,y)\ls\va\}$,
$\Bo(x,\va)=\{y\in T:d(x,y)<\va\}$. The diameter of $T$, i.e. $\sup\{d(s,t):\;s,t\in T\}$ is denoted by $D(T)$.
We define the minorizing metric
$$
\tau_{m,\varphi}(s,t):=\max\{\int^{d(s,t)}_0\varphi^{-1}(\frac{1}{m(\Bc(s,\va))})d\va,
\int^{d(s,t)}_0\varphi^{-1}(\frac{1}{m(\Bc(t,\va))})d\va \}\;\;\for\;s,t\in T.
$$
Kwapien and Rosinski \cite{Kwa} introduced these metrics to prove results
on H{\"o}lder continuity of random processes with bounded increments. However their method requires that $\varphi$ verifies (\ref{pop})
which means the exponential growth of $\varphi$. The goal of this paper is to obtain similar results, yet under relaxed conditions
imposed on $\varphi$.
\bt\label{pir1}
Let $\varphi$ and $\psi$ be Young functions (verifying $\varphi(1)=\psi(1)=1$) and for some $R>1$, $n_0\gs 1$, $n_0\in \N$
\begin{eqnarray}\label{miau}
&&\frac{\varphi(R^k)}{\varphi(R^{k+1})}\ls \frac{\varphi(R^{k-1})}{\varphi(R^{k})},\;\;\for\; k\gs 1,\;k\in \N .\\ 
\label{pap} && \sum^{\infty}_{k=0}\frac{\varphi(R^k)}{\psi(R^{k+n_0})}<\infty.
\end{eqnarray}
Let $\psi_{+}(x)=(\psi(x)-1)_{+}$ for all $x>0$. Then there exists a Borel probability measure $\nu$ on $T\times T$ and a constant $0<K<\infty$ only depending on $(\varphi,\psi)$ such that
for every continuous function $f:T\ra \R$ there holds
\be\label{1claim}
|f(s)-f(t)|\ls K|f^d|^{\nu}_{\psi_{+}}\tau_{m,\varphi}(s,t),\;\for\;s,t\in T,\;\mbox{where}\;f^d(u,v)=\frac{|f(u)-f(v)|}{d(u,v)}. 
\ee
and if $\psi\in \nabla'$, then we have
\be\label{2claim}
\sup_{s,t\in T}\psi_{+}(\frac{|f(s)-f(t)|}{Kr\tau_{m,\varphi}(s,t)})\ls \int_{T\times T}\psi_{+}(\frac{|f(u)-f(v)|}{d(u,v)})\nu(du,dv),
\ee
where $r$ is chosen such that condition (\ref{pol}) holds with $c=1$.
\et
Theorem \ref{pir1} has an application to the stochastic analysis. We say that process $X(t)$, $t\in T$ has $\varphi$-bounded increments
if it verifies
\be\label{war2}
\sup_{s,t\in T}\E\varphi(\frac{|X(s)-X(t)|}{d(s,t)})\ls 1.
\ee
\begin{colo}\label{ubu}
Suppose $(\varphi,\psi)$ verify conditions (\ref{miau}) and (\ref{pap}).
For each separable stochastic process $X(t)$, $t\in T$ which has $\psi$-bounded increments 
there holds
$$
\E\sup_{s,t\in T}\frac{|X(s)-X(t)|}{2K\tau_{m,\varphi}(s,t)}\ls 1
$$
and if $\psi\in \nabla'$ then also
$$
\E\sup_{s,t\in T}\psi(\frac{|X(s)-X(t)|}{2K\tau_{m,\varphi}(s,t)})\ls 1
$$
where $K$ is the same constant as in Theorem \ref{pir1}.
\end{colo}
\begin{dwd}
Following arguments from the proof of Theorem 2.3 in Talagrand \cite{Tal0}
it is enough to prove the result assuming that
$X(t)$, $t\in T$ has a.s. continuous samples.
Theorem \ref{pir1}, namely (\ref{1claim})  the Fubini theorem and the definition of $|\cdot|^{\nu}_{\psi_{+}}$ give
$$
\E\sup_{s,t\in T}\frac{|X(s)-X(t)|}{K\tau_{m,\varphi}(s,t)}\ls 1+\E\int_{T\times T}\psi_{+}(\frac{|X(u)-X(v)|}{d(u,v)})\nu(du,dv)\ls 2.
$$
It proves the first thesis. If $\psi\in \nabla'$, then we can apply (\ref{2claim}) instead of (\ref{1claim}) obtaining  
\ba
&& \E\sup_{s,t\in T}\psi(\frac{|X(s)-X(t)|}{K\tau_{m,\varphi}(s,t)})\ls
1+\E\sup_{s,t\in T}\psi_{+}(\frac{|X(s)-X(t)|}{K\tau_{m,\varphi}(s,t)})\ls\nonumber\\
&& \ls 1+\E\int_{T\times T}\psi_{+}(\frac{|X(u)-X(v)|}{d(u,v)})\nu(du,dv)\ls 2.\nonumber
\ea
By the convexity of $\varphi$, we derive the second claim.
\end{dwd}
\begin{rema}
Note that if $\sum^{\infty}_{k=0}\frac{\varphi(R^k)}{\varphi(R^{k+\tau})}<\infty$, for some $R>1,n_0\gs 1$ then we can take
$\psi\equiv \varphi$ in Theorem \ref{pir1}. Thus all processes which verify (\ref{war2}) (for $\varphi$) are H{\"o}lder continuous 
with respect to $\tau_{m,\varphi}(s,t)$. If $\varphi(x)\equiv x^p$ we can take $\psi(x)\equiv x^{p+\va}$, where $\va>0$ and consequently 
obtain a generalization of basic Kolmogorov result \cite{Slu}.
\end{rema}
We then prove the converse statement that minorizing metrics are optimal when considering H{\"o}lder continuity of processes with bounded increments.
\bt\label{pir2}
Assume $(\varphi,\psi)$ verify for some $R$, $n_0\gs 1$
\be\label{pep}
\sum^{\infty}_{k=0}\frac{\psi(R^k)}{\varphi(R^{k+n_0})}<\infty.
\ee
Suppose $\rho$ is a metric on $T$ such that for each separable
process $X(t)$, $t\in T$ which has $\psi$-bounded increments (verifies condition (\ref{war2}) for $\psi$), we have
$$
\P(\sup_{s,t\in T}\frac{|X(s)-X(t)|}{\rho(s,t)}< \infty)=1,
$$
then there exist a constant $K$ and a Borel probability measure $m$
(which depends on $(\varphi,\psi)$ only) such that
$\tau_{m,\varphi}(s,t)\ls K\rho(s,t)$.
\et
\begin{rema}
If $\sum^{\infty}_{k=0}\frac{\varphi(R^k)}{\varphi(R^{k+n_0})}<\infty$ then we can take $\psi=\varphi$ in Theorem \ref{pir2}.
That means there exists $m\in\ccP(T)$ such $\tau_{m,\varphi}(s,t)\ls K\rho(s,t)$ for each  
$\rho$ with respect to which all process with $\varphi$-bounded increments are H{\"o}lder continuous.
\end{rema}
\smallskip

\noindent
We also prove some generalization of Talagrand's
Theorem 4.2 \cite{Tal0}  and the author's Theorem 1 in \cite{Bed}.
\begin{theo}\label{poly1}
Assume that $\varphi$ verifies (\ref{miau}) for some $R>1$.
There exist constants $C,K$ (depending on $\varphi$ only) and a Borel probability measure $\nu$ on $T\times T$ 
such that for each continuous function $f$ on $T$ the inequality holds
\ba
&&\sup_{s,t\in T}\varphi_{+}(\frac{|f(s)-f(t)|}{C\tau_{m,\varphi}(s,t)\varphi^{-1}_{+}(\frac{\ccM(m,\varphi)}{K\tau_{m,\varphi}(s,t)})})\ls
\int_{T\times T}\varphi_{+}(\frac{|f(u)-f(v)|}{d(u,v)})\nu(du,dv),\nonumber
\ea
where $\ccM(m,\varphi):=\int_T\int^{D(T)}_0\varphi^{-1}(\frac{1}{m(\Bc(t,\va))})d\va m(dt)<\infty.$
\end{theo}
\begin{colo}
For each separable process $X(t)$, $t\in T$ which satisfies
(\ref{war2}) (for $\varphi$) there holds
$$
\E\sup_{s,t\in T}\varphi(\frac{|X(s)-X(t)|}{C\tau_{m,\varphi}(s,t)\varphi^{-1}_{+}(\frac{\ccM(m,\varphi)}{K\tau_{m,\varphi}(s,t)})})\ls 1.
$$
\end{colo}
\begin{dwd}
As in the proof of Corollary \ref{ubu} it is enough to show the result for $X(t)$, $t\in T$ with a.s. continuous samples.
Note that $\varphi(x)\ls 1+\varphi_{+}(x)$, thus due to Theorem \ref{poly1} the Fubini theorem we obtain
$$
\E\sup_{s,t\in T}\varphi(\frac{|X(s)-X(t)|}{C\tau_{m,\varphi}(s,t)\varphi^{-1}_{+}(\frac{\ccM(m,\varphi)}{K\tau_{m,\varphi}(s,t)})})\ls
1+\int_{T\times T}\E\varphi(\frac{|X(u)-X(v)|}{d(u,v)}) \nu(du,dv)\ls 2.
$$
Now by the convexity we establish the result.
\end{dwd}
In the paper we follow methods from \cite{Bed}. For a completeness we repeat from there
some of the arguments.

\section{Notation and Preliminaries}

\noindent
{\bf Young functions}
\smallskip

\bl\label{1lem}
There holds $|h|^{\mu}_{\varphi}\ls \|h\|^{\mu}_{\varphi}\ls 2|h|^{\mu}_{\varphi}$ for every $h\in \mathfrak{B}(X)$.
\el
\begin{dwd}
First note either $\int_X\varphi(\frac{|h|}{a})d\mu\ls 1$ or $\int_X\varphi(\frac{|h|}{a})d\mu> 1$ and in this case
using that $\alpha\ra\alpha\varphi(\frac{x}{\alpha})$ is decreasing we derive  
$$
\int_X \varphi(\frac{|h|}{a\int_X\varphi(\frac{|h|}{a})d\mu})d\mu\ls 
\frac{\int_X \varphi(\frac{|h|}{a})d\mu}{\int_X \varphi(\frac{|h|}{a})d\mu}=1.
$$
Consequently $|h|^{\mu}_{\varphi}\ls a+a\int_X\varphi(\frac{|h|}{a})d\mu$ for all $a>0$. That means $|h|^{\mu}_{\varphi}\ls \|h\|^{\mu}_{\varphi}$.
The last inequality follows by taking $a=|h|^{\mu}_{\varphi}$ in the definition of $\|h\|^{\mu}_{\varphi}$.
\end{dwd}
\bl\label{2lem}
Let $\varphi$ be a Young function satisfying condition (\ref{pol}) with $c=0$ and $r>0$. Then we have
$\varphi(\frac{1}{r}|h|^{\mu}_{\varphi})\ls \int_S\varphi(|h|)d\mu$ for every $h\in \mathfrak{B}(X)$.
\el
\begin{dwd}
If $\int_S\varphi(|h|)d\mu$ is either $0$ or $\infty$, then the inequality holds trivially. Suppose that $0<\int_X\varphi(|h|)d\mu<\infty$
and let us take $C>0$ so that $\varphi(C)=\int_X\varphi(|h|)d\mu$. By (\ref{pol}) property we have $\varphi(C)\varphi(\frac{x}{rC})\ls \varphi(x)$ for all $x\gs 0$ and consequently
$$
\int_X\varphi(\frac{|h|}{rC})d\mu\ls \frac{1}{\varphi(C)}\int_X\varphi(|h|)d\mu=1.
$$
Hence, we see that $\|h\|^{\mu}_{\varphi}\ls rC$ which proves the lemma.
\end{dwd}
Observe that for each Young function $\varphi$ there holds  
\be\label{wir}
\frac{x}{y}\ls \frac{\varphi(x)}{\varphi(y)},\;\;\for\;\frac{x}{y}\gs 1.
\ee
\bl\label{hit1} If $\varphi$ satisfies (\ref{miau}) then $\varphi\in\nabla'$ with $r=R^2$ and $c=1$.
\el
\begin{dwd}
By (\ref{miau}) we have
$$
\frac{\varphi(R^k)}{\varphi(R^{k+1})}\ls \frac{\varphi(R^{k-1})}{\varphi(R^{k+1})},\;\;\for\;k\gs 1,
$$
Let $i,j\gs 0$ be such that $R^{i}\ls x< R^{i+1}$ and $R^{j}\ls y<R^{j+1}$.
Clearly
$$
\frac{\varphi(R^{i+1})}{\varphi(R^{i+1}R^{j+1})}=
\frac{\varphi(R^{i+1}R^{j})}{\varphi(R^{i+1}R^{j+1})}...\frac{\varphi(R^{i+1})}{\varphi(R^{i+2})}\ls
 \frac{\varphi(R^{j})}{\varphi(R^{j+1})}...\frac{\varphi(R^{0})}{\varphi(R^{1})}
=\frac{1}{\varphi(R^{j+1})}
$$
and hence $\varphi(x)\varphi(y)\ls \varphi(R^{i+1})\varphi(R^{j+1})\ls \varphi(R^{i+1}R^{j+1})\ls \varphi(R^2xy)$.
\end{dwd}
\smallskip

\noindent
{\bf The main construction}
\smallskip

\noindent
Fix any $R>2$. For $k\gs 0$ and $x\in T$ we define $r_0(x)=D(T)$ and
\be\label{em2}
r_{k}(x):=\min\{\va\gs 0:\;\frac{1}{m(\Bc(x,\va))}\ls \varphi(R^k)\}.
\ee
Let us notice that $r_k\ls D(T)$, for $k\gs 0$.
\begin{lema}\label{moja}
The functions $r_k$ verify the Lipschitz condition with constant $1$.
\end{lema}
\begin{dwd}
Clearly $r_0$ is a constant function so it is $1$-Lipschitz. For $k>0$ and $s,t\in T$ it is
$$
\frac{1}{m(\Bc(s,r_k(t)+d(s,t))}\ls \varphi(R^k),\;\;\mbox{and}\;\;\frac{1}{m(\Bc(t,r_k(s)+d(s,t))}\ls \varphi(R^k).
$$
Hence $r_k(s)\ls r_k(t)+d(s,t)$, $r_k(t)\ls r_k(s)+d(s,t)$, thus $r_k$ is $1$-Lipschitz.
\end{dwd}
Lemma \ref{moja} gives  that $r_k\in C(T)$. 
\begin{rema}\label{hhh}
Note that if $r(x):=\lim_{k\ra\infty}r_k(x)$, we have
$r(x)=\inf\{\va\gs 0:\;m(B(x,\va))>0\}=\mathrm{ess}\inf d(x,\cdot)$ where the essential infimum is taken with respect to the probability measure $m$. In particular $r(x)=0$ if and only if $x\in \mathrm{supp}(m)$. 
\end{rema}
For each positive integer $c$ we have
\ba
&&\frac{R-1}{R}\sum_{k\gs c}r_k(x)R^k\ls
\sum_{k\gs c}r_k(x)(R^k-R^{k-1})\ls
\sum_{k\gs c}(r_k(x)-r_{k+1}(x))R^k+\nonumber\\
&&+\limsup_{k\ra\infty}r_{k+1}(x)R^{k+1}\ls \sum_{k\gs c}\int^{r_{k}(x)}_{r_{k+1}(x)}\varphi^{-1}(\frac{1}{m(\Bc(x,\va))})d\va+\nonumber\\
&&+\limsup_{k\ra\infty}\int^{r_{k+1}(x)}_{0}\varphi^{-1}(\frac{1}{m(\Bc(x,\va))})d\va=\int^{r_c(x)}_{0}\varphi^{-1}(\frac{1}{m(\Bc(x,\va))})d\va\nonumber.
\ea
Thus
\be\label{em1}
\sum_{k\gs c}r_k(x)R^k\ls
\frac{R}{R-1}\int^{r_c(x)}_{0}\varphi^{-1}(\frac{1}{m(\Bc(x,\va))})d\va.
\ee
Let us abbreviate $\Bc(x,r_k(x))$ by $\Bc_k(x)$ and $B^{\circ}(x,r_k(x))$ by $B^{\circ}_{k}(x)$ for $k>0$.
For $k=0$ we put $B^{\circ}_{0}(x)=\Bc_0(x)=T$. Due to (\ref{em2}) it is clear that
\be\label{moja1}
\frac{1}{m(\Bc_k(x))}\ls \varphi(R^k)\ls \frac{1}{m(B^{\circ}_k(x))},\;\;\for\; k\gs 0.
\ee
For each $k\gs 0$ we define the linear operator $S_k:\mathfrak{B}_b(T)\ra \mathfrak{B}_b(T)$ by the formula
$$
S_k f(x):=\pint_{\Bc_k(x)}f(u)m(du)=\frac{1}{m( \Bc_k(x) ) }\int_{\Bc_k(x)}f(u)m(du).
$$
If $f,g\in\mathfrak{B}_b(T)$, $k\gs 0$, then we easily check that:
\begin{itemize}
\item[(i)] $S_k 1=1$;
\item[(ii)] if $f\ls g$ then $S_k f\ls S_k g$ and hence $|S_k f|\ls S_k |f|$;
\item[(iii)] if $f\in C(T)$ and $\lim_{k\ra\infty}r_k(x)=0$, then $\lim_{k\ra\infty}S_k f(x)=f(x)$.
\end{itemize}
Fix $l\gs 0$. There exists unique $m^l_{x,k}\in \ccP(T)$ such that
for each $f\in \mathfrak{B}_b(T)$ we have
\be\label{nim}
S_l S_{l-1}... S_{k}f(x)=\int_T f(u) m^l_{x,k}(du),\;\;\for\; 0\ls k\ls l.
\ee
Let us define
$$
r^l_k:=\sum^{l}_{i=k}2^{i-k}r_i,\;\;\Bc^l_k(x):=\Bc(x,r^l_k(x)),\;\;\for\;k\ls l.
$$
\begin{lema}\label{be1}
For each $u\in \Bc^l_{k+1}(x)$ $0\ls k<l$ we have $\Bc_k(u)\subset \Bc^l_k(x)$ and
$$
r_{k}(u)\ls r_k(x)+r^l_{k+1}(x)\ls r^l_k(x).
$$
\end{lema}
\begin{dwd}
Fix $u\in \Bc^l_{k+1}(x)$. Since $r_k$ are $1$-Lipschitz, we get
$$
r_k(u)\ls r_k(x)+d(x,u)\ls r_k(x)+r^l_{k+1}(x)\ls r^l_{k}(x).
$$
Clearly $r_k(u)\ls r_k(x)+r^l_{k+1}(x)$. Furthermore $d(x,u)\ls r^l_{k+1}(x)$, thus
$$
\Bc(u,r_k(u))\subset \Bc(u,r_k(x)+r^l_{k+1}(x))\subset \Bc(x,r_k(x)+2r^l_{k+1}(x))=\Bc(x,r^l_k(x))
$$
and by the definition $\Bc_k(u)=\Bc(u,r_k(u))$, $\Bc^l_k(x)=\Bc(x,r^l_k(x))$.
\end{dwd}
\bl\label{dwa}
For all $0\ls k\ls l$
we have $m^l_{x,k}(\Bc^l_k(x))=1$ i.e. $\mathrm{supp}(m^l_{x,k})\subset \Bc^l_k(x)$.
\el
\begin{dwd}
We prove Lemma \ref{dwa} by the reverse induction on $k$.
Clearly $\mathrm{supp}(m^l_{x,l})=\Bc(x,r_l(x))=\Bc^l_l(x)$.
Suppose that for some $k<l$ we have $\mathrm{supp}(m^l_{x,k+1})\subset \Bc^l_{k+1}(x)$, then
the definition gives
$$
\int_T f(u)m^l_{x,k}(du)=\int_T\pint_{\Bc_k(u)}f(v)m(dv)m^l_{x,k+1}(du),\;\;\for\;f\in \mathfrak{B}_b(T).
$$
Due to Lemma \ref{be1} we have
$\Bc_k(u)\subset \Bc^l_k(x)$, for $u\in \Bc^l_{k+1}(x)$. It ends the proof.
\end{dwd}
\begin{colo}\label{piss}
For each $f\in \mathfrak{B}_b(T)$, and $k\ls l$ the inequality holds
$$
S_l S_{l-1}... S_{k}|f|(x)=\int_T |f(u)|m^l_{x,k}(du)\ls \varphi(R^k)\int_{\Bc^l_{k}(x)}|f(u)|m(du).
$$
\end{colo}
\begin{dwd}
If $k=l$ the inequality is obvious. If $k<l$, using Lemma \ref{dwa}, and (\ref{moja1}) we obtain
\ba
&&S_l S_{l-1}... S_{k}|f|(x)=\int_T\pint_{\Bc_k(u)}|f(v)|m(dv)m^l_{x,k+1}(du)\ls\nonumber\\
&&\ls \varphi(R^k)\int_T\int_{\Bc^l_k(x)}|f(v)|m(dv)m^{l}_{x,k+1}(du)=\varphi(R^k)\int_{\Bc^l_k(x)}|f(v)|m(dv).\nonumber
\ea
\end{dwd}
Let us notice that for a positive integer $c$ with $0\ls c<l$ we have
\be
\sum^{l-1}_{k=c} r^l_k R^{k}=
\sum^{l-1}_{k=c}\sum^{l}_{i=k} (\frac{2}{R})^{i-k} r_i R^i\ls \sum^\infty_{j=0}(\frac{2}{R})^j\sum^{l}_{i=c}r_i R^i=
\frac{R}{R-2}\sum^{\infty}_{i=c}r_iR^i.\nonumber
\ee
Together with (\ref{em1}) it gives
\be\label{szesc}
\sum^{l-1}_{k=c} r^{l}_k(x)R^{k}\ls \frac{R^2}{(R-1)(R-2)}\int^{r_c(x)}_0\varphi^{-1}(\frac{1}{m(\Bc(x,\va))})d\va.
\ee

\section{Proof of Theorem 1}

\begin{dwd}
We may assume that (\ref{miau}) and (\ref{pap}) hold with $R>5$ (note that if (\ref{miau}) and (\ref{pap}) hold for some $R$ then they hold also
for $R^l$, where $l\in \N$). Fix $s,t\in T$, without losing the generality we may assume also $\tau_{m,\varphi}(s,t)<\infty$, which implies that $\lim_{k\ra \infty}r_k(x)=0$, for $x=s,t$. If $d(s,t)<D(T)$ then there exist positive integers $a,b$ such that
$$
r_{a}(s)\ls  d(s,t)<  r_{a-1}(s),\;\; r_{b}(t)\ls  d(s,t)< r_{b-1}(t),
$$
and we can define $c:=\max\{a,b\}$. If $d(s,t)=D(T)=r_0$, we put $c:=0$. For a fixed $l>c$ let us denote 
$$
\tau_{x}:=\max\{k\gs 1:\;\Bc^l_k(s)\cup \Bc^l_k(t)\subset B^{\circ}_{k-1}(u),\;\;\For\;u\in \Bc^l_{k}(x)\},
\;\;x=s,t.
$$
and $\tau:=\min\{\tau_s,\tau_t\}$. Observe that $B^{\circ}_{0}(u)=T$, for all $u\in T$ so $\tau_x$ is well defined and clearly 
$1\ls \tau\ls c$. For simplicity we put also $r^l_{k}(s,t):=r^l_k(s)+r^l_k(t)$ and $d_k(s,t):=\min\{r^l_{k}(s,t)+d(s,t),D(T)\}$.
Note that 
\be\label{palinka}
d_{\tau}(s,t)\ls r_{\tau-1}(u),\;\; \For\; u\in \Bc^l_{\tau}(x)\;\;\mbox{if}\;\;\tau=\tau_x.
\ee
\bl\label{latex0}
The inequality holds 
$$
d_{\tau}(s,t)R^{\tau}+\sum^{c}_{k=\tau}R^k r^l_k(s,t) \ls \frac{R}{R-5} R^{c}(\frac{3}{2}d(s,t)+2r^l_c(s,t)).
$$
\el
\begin{dwd}
Let $\tau\ls k < c$ be given and let $x$ be either $s$ or $t$. There exist $u_x\in \Bc^l_{k+1}(x)$, $x=s,t$
such that $r_{k}(u_x)\ls d_k(s,t)$.
Indeed, otherwise
$$
\Bc^l_{k+1}(s)\cup \Bc^l_{k+1}(t)\subset \Bc(u,d_{k+1}(s,t))\subset B^\circ_{k}(u)\;\;\For\;
u\in \Bc^l_{k+1}(t)\cup \Bc^l_k(s)
$$
which is impossible due to the definition of $\tau$.
\smallskip

\noindent
By Lemma \ref{moja} functions $r_k$ are 1-Lipschitz, therefore
$$
r_{k}(x)\ls r_{k}(u_x)+r^l_{k+1}(x)\ls d_{k+1}(s,t)+r^l_{k+1}(x),\;\;x=s,t.
$$
Since $r^l_{k}=r_{k}+2r^l_{k+1}$, we obtain $r^{l}_{k}(x)\ls d_{k+1}(s,t)+3r^l_{k+1}(x)$. Consequently
$$
r^l_{k}(s,t)\ls 2d_{k+1}(s,t)+3r^l_{k+1}(s,t)= 2d(s,t)+5r^l_{k+1}(s,t).
$$
Iterating this inequality, we obtain the following result 
\be\label{ucl}
r^l_k(s,t)\ls 2d(s,t)\ls \sum^{c-k-1}_{i=0}5^i+5^{c-k}r^l_c(s,t)=\frac{d(s,t)}{2}(5^{c-k}-1)+5^{c-k}r^l_c(s,t)
\ee
for all $\tau\ls k\ls c$ (observe that inequality holds trivially for $k=c$). Hence, we have
$$
\sum^c_{k=\tau}r^l_k(s,t)\ls (\frac{d(s,t)}{2}+r^l_{c}(s,t))\sum^c_{k=\tau}R^k5^{c-k}\ls \frac{R}{R-5}R^c(\frac{d(s,t)}{2}+r^l_{c}(s,t))
$$
and by (\ref{ucl}) we have (recall that $R>5$)
\ba
& d_{\tau}(s,t)R^{\tau}\ls R^{\tau}(d(s,t)+r^l_{\tau}(s,t))\ls d(s,t)(1+\frac{1}{2}(5^{c-\tau}-1))R^{\tau}+5^{c-\tau}R^{\tau}r^l_c(s,t)\ls\nonumber\\
& \ls 5^{c-\tau}R^{\tau}(d(s,t)+r^l_c(s,t))\ls R^c(d(s,t)+r^l_c(s,t)).
\ea
Since $\frac{R}{R-5}>1$, we obtain the inequality.
\end{dwd}
We remind that $f^d(u,v)=\frac{|f(u)-f(v)|}{d(u,v)}$.  For simplicity we denote 
$$
F_k:=\{(u,v)\in T\times T:\;\;f^d(u,v)\gs R^{k}\},\;\;k\gs 0.
$$
\bl\label{letmat1}
If $\varphi$ satisfies (\ref{miau}), then for each positive integer $n$ and $f\in C(T)$ there holds
\ba
&&|S_l f(s)-S_l f(t)|\ls d_{\tau}(s,t)R^{\tau+n}+\sum_{x\in\{s,t\}}\sum^{l-1}_{k=\tau}r^l_k(x)R^{k+n}+\nonumber\\
&&+\sum_{x\in\{s,t\}}\sum^{l-1}_{k=\tau}
\varphi(R^{k+1})\int_{\Bc^l_{k+1}(x)}r_k(u)\pint_{\Bc_k(u)}f^d(u,v)1_{F_{k+n}}m(dv)m(du))+\nonumber\\
&&+d_{\tau}(s,t)\varphi(R^{\tau+1})
\int_{\Bc^l_{\tau}(y)}\pint_{B^{\circ}_{\tau-1}(u)}f^d(u,v)1_{F_{\tau+n}}m(dv)m(du)),\nonumber
\ea
where $y=t$ if $\tau=\tau_t$ and $y=s$ if $\tau\neq \tau_t$.
\el
\begin{dwd}
Fix $f\in C(T)$. Without losing the generality generality we can assume that $\tau=\tau_t$.  Clearly
\ba\label{rer0}
&& S_{l}f(s)-S_lf(t)= \sum^{l-1}_{k=\tau}S_l...S_{k+1}(\mathrm{Id}-S_{k})f(s)-\nonumber\\
&& -\sum^{l-1}_{k=\tau}S_l...S_{k+1}(\mathrm{Id}-S_{k})f(t)+(S_l...S_{\tau}f(s)-S_l...S_{\tau}f(t)).
\ea
We have also
\be\label{kaczor}
|S_l...S_{k+1}(\mathrm{Id}-S_{k})f(x)|\ls \int_T |(\mathrm{Id}-S_{k})f(u)| m^l_{x,k+1}(du),
\ee
Since $f^d(u,v)\ls R^{k+n}+f^d(u,v)1_{F_{k+n}}$,
we obtain 
\ba
&& |(\mathrm{Id}-S_{k})f(u)|\ls \pint_{\Bc_k(u)}|f(u)-f(v)|m(dv)\ls r_k(u)\pint_{\Bc_k(u)}f^d(u,v)m(dv)\ls\nonumber\\ 
&& \ls r_k(u)R^{k+n}+r_k(u)\pint_{\Bc_k(u)}f^d(u,v)1_{F_{k+n}}m(dv),\;\;\For\;u\in T.\nonumber
\ea
By Lemma \ref{be1}, $r_k(u)\ls r^l_k(x)$, whenever $u\in \Bc^l_{k+1}(x)$.
This, (\ref{kaczor}) and Corollary \ref{piss} imply that
\ba\label{rer1}
&&|S_l...S_{k+1}(\mathrm{Id}-S_{k})f(x)|\ls
\int_T |(\mathrm{Id}-S_{k})f(u)| m^l_{x,k+1}(du)\ls r^l_k(x)R^{k+n}+ \nonumber\\
&& +\int_T r_k(u)\pint_{\Bc_{k}(u)}f^d(u,v)1_{F_{k+n}}m(dv)m^l_{x,k+1}(du)\ls r^l_k(x)R^{k+n}+
\nonumber\\
&&+\varphi(R^{k+1})\int_{\Bc^l_{k+1}(x)}r_k(u)\pint_{\Bc_{k}(u)}f^d(u,v)1_{F_{k+n}}m(dv)m^l_{x,k+1}(du).
\ea
To bound the last part in (\ref{rer0}) let us observe that
\be\label{vava0}
|S_l...S_{\tau}f(s)-S_l...S_{\tau}f(t)|\ls \int_T\int_T |f(u)-S_{\tau}f(w)|m^l_{s,\tau+1}(dw)m^l_{t,\tau}(du).
\ee
By Lemma \ref{dwa}
$\mathrm{supp}(m^l_{x,k})\subset \Bc^l_{k}(x)$,
$x\in T$.
If $w\in \Bc^l_{\tau+1}(s)$ and $u\in \Bc^l_{\tau}(t)$, then
$$
|f(u)-S_{\tau}f(w)|\ls \pint_{\Bc_{\tau}(w)}|f(u)-f(v)|m(dv).
$$
Lemma \ref{be1} implies that
$\Bc_{\tau}(w)\subset \Bc^l_{\tau}(s)$. Hence for each $u\in \Bc^l_{\tau}(t)$, $v\in \Bc_{\tau}(w)$
\be\label{der2}
d(u,v)\ls \min\{d(u,t)+d(t,s)+d(s,v),D(T)\}\ls d_{\tau}(s,t).
\ee
Applying (\ref{der2}) and $f^d(u,v)\ls R^{\tau+n}+f^d(u,v)1_{F_{\tau+n}}$ we obtain 
\ba\label{vava1}
&& |f(u)-S_{\tau}f(w)|\ls d_{\tau}(s,t)
\pint_{\Bc_{\tau}(w)}f^d(u,v)m(dv)\ls \nonumber\\
&&\ls d_{\tau}(s,t)(R^{\tau+n}+\pint_{\Bc_{\tau}(w)}f^d(u,v)1_{F_{\tau+n}}m(dv)).
\ea
Since $\tau=\tau_t$ we have
$\Bc_{\tau}(w)\subset B^l_{\tau}(s)\subset B^{\circ}_{\tau-1}(u)$ for all $w\in \Bc^l_{\tau+1}(t)$.
Together with (\ref{moja1}) it implies 
\ba\label{pen1}
&&\pint_{\Bc_{\tau}(w)}f^d(u,v)1_{F_{\tau+n}}m(dv)
\ls\varphi(R^{\tau})
\int_{\Bc_{\tau}(w)}f^d(u,v)1_{F_{\tau+n}}m(dv)\ls\nonumber\\
&&\ls\frac{\varphi(R^{\tau})}{\varphi(R^{\tau-1})}
\pint_{B^{\circ}_{\tau-1}(u)}f^d(u,v)1_{F_{\tau+n}}m(dv).
\ea
The condition (\ref{miau}) gives
$\frac{\varphi(R^{\tau})}{\varphi(R^{\tau-1})}\ls \frac{\varphi(R^{\tau+1})}{\varphi(R^{\tau})}$.
Hence, due to (\ref{vava1}) and (\ref{pen1}) we obtain
\be\label{vava2}
|f(u)-S_{\tau}f(w)|\ls d_{\tau}(s,t)(R^{\tau+n}+
\frac{\varphi(R^{\tau+1})}{\varphi(R^\tau)}
\pint_{B^{\circ}_{\tau-1}(u)}f^d(u,v)1_{F_{\tau+n}}m(dv)).
\ee
Inequalities (\ref{vava0}), (\ref{vava2}) and Corollary \ref{piss} imply
\ba\label{rer3}
&& |S_l...S_{\tau}f(s)-S_l...S_{\tau}f(t)|\ls \nonumber\\
&&\ls d_{\tau}(s,t)(R^{\tau+n}+\frac{\varphi(R^{\tau+1})}{\varphi(R^\tau)}\pint_{B^{\circ}_{\tau-1}(u)}f^d(u,v)1_{F_{\tau+n}}m(dv)m^l_{t,\tau}(du))\ls \nonumber\\
&&\ls d_{\tau}(s,t)(R^{\tau+n}+\varphi(R^{\tau+1})\int_{\Bc^l_\tau(t)}\pint_{B^{\circ}_{\tau-1}(u)}
f^d(u,v)1_{F_{\tau+n}}m(dv)m(du)).
\ea
Note that (\ref{rer1}) and (\ref{rer3}) give the result
\end{dwd}
\bl\label{letmat2}
If $A=\frac{4R^3}{(R-1)(R-2)(R-5)}+\frac{3R^2}{2(R-5)}$, then we have
$$
d_{\tau}(s,t)R^{\tau}+\sum_{x\in\{s,t\}}
\sum^{l-1}_{k=\tau}r^l_k(x)R^{k}\ls A\tau_{m,\varphi}(s,t).
$$
\el
\begin{dwd}
Lemma \ref{latex0} gives
\ba
&&d_{\tau}(s,t)R^{\tau}+\sum_{x\in\{s,t\}}\sum^{l-1}_{k=\tau}r^l_k(x)R^{k}
=\sum^c_{k=\tau}r^l_k(s,t)R^{k}+\sum^{l-1}_{k=c+1}r^l_{k}(s,t)R^k\ls\nonumber\\
&&\ls \frac{R}{R-5}(\frac{3}{2}d(s,t)+2\sum^{l-1}_{k=c}r^l_k(s,t)R^k).
\nonumber
\ea
Clearly $r_c(x)\ls d(s,t)$, $x\in\{s,t\}$, thus by (\ref{szesc}) we obtain 
$$ 
2\sum^{l-1}_{k=c}(r^l_k(s)+r^l_k(t))R^{k}\ls \frac{4R^2}{(R-1)(R-2)}\max_{x\in\{s,t\}}\int^{d(s,t)}_0\varphi^{-1}(\frac{1}{m(\Bc(x,\va))})d\va.
$$
Since $d(s,t)< \max\{r_{c-1}(s),r_{c-1}(t)\}$ if $c>0$ and $d(s,t)=D(T)$ if $c=0$,
we have
$$
R^{c-1}\ls \max_{x\in\{s,t\}}\varphi^{-1}(\frac{1}{m(\Bc(x,d(s,t))}).
$$
It follows that
$$
d(s,t)R^{c}\ls R\max_{x\in\{s,t\}}\int^{d(s,t)}_0\varphi^{-1}(\frac{1}{m(\Bc(x,\va))})d\va.
$$
Hence, due to the definition of $\tau_{m,\varphi}(s,t)$ we deduce
$$
d_{\tau}(s,t)R^{\tau}+\sum_{x\in\{s,t\}}\sum^{l-1}_{k=\tau}r^l_k(x)R^{k}\ls A\tau_{m,\varphi}(s,t).
$$
\end{dwd}
Lemma \ref{be1} implies
$r_k(u)\ls r^l_k(x)$, for $u\in B^l_k(x)$. This observation together with Lemma \ref{letmat1} (with $n=n_0+1$) yields 
\ba
&&|S_l f(s)-S_l f(t)|\ls d_{\tau}(s,t)R^{\tau+n_0+1}+\sum_{x\in\{s,t\}}\sum^{l-1}_{k=\tau}r^l_k(x)R^{k+n_0+1}+\nonumber\\
&&+\sum_{x\in\{s,t\}}\sum^{l-1}_{k=\tau}r^l_k(x)R^{k+n_0+1}\varphi(R^{k+1})\int_{\Bc^l_{k+1}(x)}\pint_{\Bc_k(u)}\frac{f^d(u,v)}{R^{k+n}}1_{F_{k+n}}m(dv)m(du)+\nonumber\\
&&+d_{\tau}(s,t) R^{\tau+n_0+1}\varphi(R^{\tau+1})
\int_{\Bc^l_{\tau}(y)}\pint_{B^{\circ}_{\tau-1}(u)}\frac{f^d(u,v)}{R^{\tau+n_0+1}}1_{F_{\tau+n_0+1}}m(dv)m(du).\nonumber
\ea
By Lemma \ref{letmat2} we obtain
\ba\label{main}
&&|S_l f(s)-S_l f(t)|\ls AR^{n_0+1}\tau_{m,\varphi}(s,t)(1+\nonumber\\
&&+\sum_{x\in\{s,t\}}\sum^{\infty}_{k=1}
\varphi(R^{k+1})\int_{T}\pint_{\Bc_k(u)}\frac{f^d(u,v)}{R^{k+n_0+1}}1_{F_{k+n_0+1}}m(dv)m(du)+\nonumber\\
&&+\sum^{\infty}_{k=1}\varphi(R^{k+1})
\int_{T}\pint_{B^{\circ}_{k-1}(u)}\frac{f^d(u,v)}{R^{k+n_0+1}}1_{F_{k+n_0+1}}m(dv)m(du)).
\ea
For each $k\gs 0$ applying (\ref{wir}) (for $\psi$) we have
\be\label{mrau}
\frac{f^d(u,v)}{R^{k}}1_{F_{k}}\ls\frac{1}{\psi_{+}(R^{k})}\psi_{+}(f^d(u,v))\ls\frac{1}{\psi_{+}(R^{k})}\psi_{+}(f^d(u,v)).
\ee
The right hand side of (\ref{main}) does not depend on $l$, furthermore the property (iii) of $S_l$ gives that $\lim_{l\ra\infty}S_lf(x)=f(x)$, for $x\in\{s,t\}$.
Hence combining (\ref{mrau}) and (\ref{main}) we obtain
\ba\label{fer1}
&&\frac{|f(s)-f(t)|}{AR^{n_0+1}\tau_{m,\varphi}(s,t)}\ls
1+2\sum^{\infty}_{k=1}\frac{\varphi(R^{k+1})}{\psi_{+}(R^{k+n_0+1})}
\int_T\pint_{\Bc_{k}(u)}\psi_{+}(f^d(u,v))m(dv)m(du)+\nonumber\\
&&+\sum^{\infty}_{k=1}\frac{\varphi(R^{k+1})}{\psi_{+}(R^{k+n_0+1})}
\int_T\pint_{B^\circ_{k-1}(u)}\psi(f^d(u,v))1_{F_0}m(dv)m(du).
\ea
It remains to construct a suitable $\nu\in \ccP(T\times T)$. For each $g\in C(T\times T)$ we put
\ba
&&\nu(g):=\frac{1}{B}\sum^{\infty}_{k=1}\frac{\varphi(R^{k+1})}{\psi_{+}(R^{k+n_0+1})}
(2\int_T\pint_{\Bc_{k}(u)}g(u,v)m(dv)m(du)+\nonumber\\
&&+\int_T\pint_{B^\circ_{k-1}(u)}g(u,v)m(dv)m(du)),\nonumber
\ea
where $B$ is such that $\nu(1)=1$. This constant exists due to (\ref{pap}), indeed
\ba
&& B=3\sum^{\infty}_{k=1}\frac{\varphi(R^{k+1})}{\psi_{+}(R^{k+n_0+1})}= 3\sum^{\infty}_{k=1}\frac{\varphi(R^k)}{\psi(R^{k+n_0+1})-1}\ls \nonumber\\
&&\ls \frac{3}{1-R^{-n_0-1}}\sum^{\infty}_{k=1}\frac{\varphi(R^k)}{\psi(R^{k+n_0+1})}<\infty,\nonumber
\ea
where we have used that $\psi(x)\ls \psi_{+}(x)+1$ and $\psi(R^{k+n_0+1})-1\gs (1-R^{-n_0-1})\psi(R^{k+n_0+1})$ (by convexity).
Plugging $\nu$ in (\ref{fer1}) and then using homogeneity, we see
\be\label{pot1}
\frac{|f(s)-f(t)|}{ABR^{n_0+1}|f^d|^{\nu}_{\psi_{+}}\tau_{m,\varphi}(s,t)}\ls
1+2\int_{T\times T}\psi_{+}(\frac{f^d(u,v)}{|f^d|^{\nu}_{\psi_{+}}})\nu(du,dv)\ls 3.
\ee
Thus we obtain (\ref{1claim}) with $K=3ABR^{n_0+1}$.
Suppose now that $\psi(x)\psi(y)\ls \psi(rxy)$ for all $x,y\gs 1$. Since $\psi(x)\gs \psi(1)=1$ for all $x\gs 1$, we have 
$\psi_{+}(x)\psi_{+}(y)\ls \psi_{+}(rxy)$ for all $x,y\gs 0$ and so we see that (\ref{2claim}) follows from (\ref{1claim})
and Lemma \ref{2lem}.
\end{dwd}

\section{Proof of Theorem  2}

\begin{dwd}
We give a proof which modifies the idea from the paper \cite{Kwa}.
In the same way as Theorem 2.3 in \cite{Tal0} it can be proved that the existence of metric $\rho$ on $T\times T$ such that for each separable process $X(t)$, $t\in T$ which satisfies (\ref{war2}) (for $\psi$) there holds
$$
\P(\sup_{s,t\in T}\frac{|X(s)-X(t)|}{\rho(s,t)}<\infty)=1,
$$
implies the existence of a constant $K_0$ and a continuous positive functional $\Lambda$ on $C_b(T\times T\backslash \triangle)$
(where $\triangle:=\{(t,t):\;t\in T\}$) with $\Lambda(1)=1$ such that for each $f\in C(T)$
\be\label{saser}
\sup_{s,t\in T}\frac{|f(s)-f(t)|}{K_0\rho(s,t)}\ls 1+\Lambda(\psi(f^d)),
\ee
where $f^d(u,v)=\frac{|f(u)-f(v)|}{d(u,v)}$.
We define measure $m\in \ccP(T)$ by the requirement
\be\label{ppp}
\int_T g(t)m(dt)=\Lambda(\frac{g(u)+g(v)}{2}),\;\;\for\;g\in C(T).
\ee
Fix $s,t\in T$ and $l\in \N$. Let us denote
$$
h_l(\va):=
\left\{\begin{array}{lll}
R^{-n_0}  & r_1(t)\ls \va \ls r_0(t) \\
R^{k-n_0} & r_{k+1}(t)\ls \va< r_k(t),\;\; 0<k\ls l\\
0       & 0\ls \va\ls r_{l+1}(t),
\end{array}\right.
$$
where $r_k(x)=\min\{\va:\;\frac{1}{m(x,\va)}\ls\varphi(R^k)\}$, for $k\gs 0$ as in our main construction.
Observe that $h_l$, $l\gs 1$ is an increasing family of functions, so $h:=\lim_{l\ra \infty}h_l$ is well defined.
We denote $f_l(x):=\int^{d(t,x)}_0 h_l(\va)d\va$ and observe that
$$\frac{|f_l(u)-f_l(v)|}{d(u,v)}\ls \frac{1}{|d(t,u)-d(t,v)|}
|\int^{d(t,u)}_{d(t,v)}h_l(\va)d\va|=
|\pint^{d(t,u)}_{d(t,v)}h_l(\va)d\va|.
$$
The Jensen's inequality gives
$$
\psi(\frac{|f_l(u)-f_l(v)|}{d(u,v)})\ls |\pint^{d(t,u)}_{d(t,v)}\psi(h_l(\va))d\va|\ls
\psi(h_l(d(t,u)))+\psi(h_l(d(t,v))),
$$
thus by (\ref{ppp}) we have
\be\label{mall}
\Lambda(\psi(f^d_l))\ls 2\int_{T}\psi(h_l(d(t,u)))m(du).
\ee
Using the definition of $h_l$ and (\ref{moja1}) we obtain
\be\label{spi}
\int_T \psi(h_l(d(t,u)))m(du)=\sum^l_{k=0}\psi(R^{k-n_0})m(B^{\circ}_k(t) \backslash B^{\circ}_{k+1}(t))
\ls \sum^{l}_{k=0}\frac{\psi(R^{k-n_0})}{\varphi(R^k)}.
\ee
Applying (\ref{pep}) we derive $D:=\sum^{\infty}_{k=0}\frac{\psi(R^{k-n_0})}{\varphi(R^k)}<\infty$.
Consequently (\ref{saser}), (\ref{mall}), (\ref{spi}) yield
$$
\frac{\int^{d(s,t)}_0 h_l(\va)d\va}{K_0\rho(s,t)}\ls 1+\Lambda(\psi(f^d_l))
\ls 1+2D.
$$
The right hand side does not depend on $l$, so 
\be\label{synek}
\frac{\int^{d(s,t)}_0 h(\va)d\va}{K_0\rho(s,t)}\ls 1+2D.
\ee
The definition of $h$ gives
$$
\varphi^{-1}(\frac{1}{m(B(t,\va))})\ls R^{k+1}=R^{n_0+1} h(\va),\;\;\for\;r_{k+1}(t)\ls \va< r_k(t),
$$
thus for $\delta\in [r_{k+1}(t),r_k(t))$, $k\in\N$
$$
R^{-n_0-1}\int^{\delta}_{r_{k+1}(t)}
\varphi^{-1}(\frac{1}{m(B(t,\va))})d\va\ls\int^{\delta}_{r_{k+1}(t)}h(\va)d\va
$$
and hence due to (\ref{synek}) we obtain
$$
\int^{d(s,t)}_0 \varphi^{-1}(\frac{1}{m(B(t,\va))})d\va\ls K\rho(s,t),
$$
where $K=(1+2D)R^{n_0+1}K_0$. Similarly
$$
\int^{d(s,t)}_0\varphi^{-1}(\frac{1}{m(B(s,\va))})d\va\ls K\rho(s,t),
$$
which means $\tau_{m,\varphi}(s,t)\ls K\rho(s,t)$.
\end{dwd}

\section{Proof of Theorem 3}

\begin{dwz}
Fix $R>5$, $s,t\in T$ and $f\in C(T)$. We can assume that $\tau_{m,\varphi}(s,t)<\infty$
which implies $\lim_{k\ra\infty}r_k(x)=0$ for $x=s,t$.
By Lemma \ref{letmat1} (with $n=1$) and (\ref{palinka}) we have
\ba
&&|S_l f(s)-S_l f(t)|\ls d_{\tau}(s,t)R^{\tau+1}+\sum_{x\in\{s,t\}}\sum^{l-1}_{k=\tau}r^l_k(x)R^{k+1}+\nonumber\\
&&+\sum_{x\in\{s,t\}}\sum^{l-1}_{k=\tau}
\varphi(R^{k+1})\int_{\Bc^l_{k+1}(x)}r_k(u)\pint_{\Bc_k(u)}f^d(u,v)1_{F_{k+1}}m(dv)m(du)+\nonumber\\
&&+\varphi(R^{\tau+1})
\int_{\Bc^l_{\tau}(y)}r_{\tau-1}(u)\pint_{B^{\circ}_{\tau-1}(u)}f^d(u,v)1_{F_{\tau+1}}m(dv)m(du),\nonumber
\ea
where $y=t$ if $\tau=\tau_t$ and $y=s$ if $\tau\neq \tau_t$. By Lemma \ref{letmat2} we obtain
\ba\label{main10}
&&|S_l f(s)-S_l f(t)|\ls AR\tau_{m,\varphi}(s,t)+\nonumber\\
&&+\sum_{x\in\{s,t\}}\sum^{\infty}_{k=1}
\varphi(R^{k+1})\int_{T}r_k(u)R^{k+1}\pint_{\Bc_k(u)}\frac{f^d(u,v)}{R^{k+1}}1_{F_{k+1}}m(dv)m(du)+\nonumber\\
&&+\sum^{\infty}_{k=1}\varphi(R^{k+1})
\int_{T}r_{k-1}(u)R^{k+1}\pint_{B^{\circ}_{k-1}(u)}\frac{f^d(u,v)}{R^{k+1}}1_{F_{k+1}}m(dv)m(du).
\ea
The condition (\ref{wir}) gives that for each $k\gs 0$
\be\label{mrau10}
\frac{f^d(u,v)}{R^{k}}1_{F_{k}}\ls\frac{1}{\varphi_{+}(R^k)}\varphi(f^d(u,v))1_{F_k}
\ls\frac{1}{\varphi_{+}(R^{k})}\varphi_{+}(f^d(u,v)).
\ee
The right hand side of (\ref{main10}) does not depend on $l$ thus we can take 
the limit on left-hand side which is $\lim_{l\ra\infty}S_lf(x)=f(x)$, for all $x\in T$
(by property (iii) of $S_l$). Observe also that by the convexity $\varphi_{+}(R^{k+1})-1\gs (1-R^{-1})\varphi(R^{k+1})$.
Consequently due to (\ref{main10}) and (\ref{mrau10}) we obtain
\ba\label{ferdek}
&&\frac{|f(s)-f(t)|}{AR\tau_{m,\varphi}(s,t)}\ls
1+\frac{1}{1-R^{-1}}(2\sum^{\infty}_{k=1}\int_T r_k(u)R^{k+1}
\pint_{\Bc_{k}(u)}\varphi_{+}(f^d(u,v))m(dv)m(du)+\nonumber\\
&&+\sum^{\infty}_{k=1}\int_T r_{k-1}(u)R^{k+1}\pint_{B^{\circ}_{k-1}(u)}\varphi_{+}(f^d(u,v))m(dv)m(du)).
\ea
To construct a probability measure $\nu\in \ccP(T\times T)$ we put for each $g\in C(T\times T)$
\ba
&&\nu(g):=\frac{1}{M(1-R^{-1})}\sum^{\infty}_{k=1}(
2\int_T r_{k}(u)R^{k+1}\pint_{\Bc_{k}(u)}g(u,v)m(dv)m(du)+\nonumber\\
&&+\int_T r_{k-1}(u)R^{k+1}\pint_{B^{\circ}_{k-1}(u)}g(u,v)m(dv)m(du),\nonumber
\ea
where $M$ is such that $\nu(1)=1$. Applying (\ref{em1}) and the definition $\ccM(m,\varphi)$ we get
\ba
&& 1=\frac{1}{M(1-R^{-1})}\sum^{\infty}_{k=1}(2\int_T r_{k}(u)R^{k+1}m(du)+\int_T r_{k-1}(u)R^{k+1}m(du))\ls\nonumber\\
&&\ls \frac{3}{M(1-R^{-1})}\sum^{\infty}_{k=0}\int_T r_{k}(u)R^{k+2}m(du)\ls \frac{3R^4}{M(R-1)^2}\ccM(m,\varphi).
\nonumber
\ea
Hence $M\ls B\ccM(m,\varphi)$, where $B=\frac{3R^4}{(R-1)^2}$. Plugging $\nu$ into (\ref{ferdek}) we obtain
$$
|f(s)-f(t)|\ls AR\tau_{m,\varphi}(s,t)
+B\ccM(m,\varphi)\int_{T\times T}\varphi_{+}(f^d(u,v)) \nu(du,dv).
$$
By homogeneity we obtain for all $a>0$
\be\label{1pot1}
\frac{|f(s)-f(t)|}{aR^2|f^d|^{\nu}_{\varphi_{+}}}\ls AR\tau_{m,\varphi}(s,t)
+B\ccM(m,\varphi)\int_{T\times T}\varphi_{+}(\frac{f^d(u,v)}{aR^2|f^d|^{\nu}_{\varphi_{+}}}) \nu(du,dv).
\ee
Due to Lemma \ref{hit1} we know that $\varphi\in\nabla'$ with $r=R^2$ and $c=1$, thus $\varphi_{+}\in \nabla'$ with $c=0$ and $r=R^2$.
Consequently by (\ref{pol}) we get
$$
\varphi_{+}(a)\int_{T\times T}\varphi_{+}(\frac{f^d(u,v)}{aR^2|f^d|^{\nu}_{\varphi_{+}}}) \nu(du,dv)\ls
\int_{T\times T}\varphi_{+}(\frac{f^d(u,v)}{|f^d|^{\nu}_{\varphi_{+}}})\nu(du,dv)=1. 
$$
Using the above inequality in (\ref{1pot1}) we obtain
$$
\frac{|f(s)-f(t)|}{aR^2|f^d|^{\nu}_{\varphi_{+}}}\ls AR\tau_{m,\varphi}(s,t)+\frac{B\ccM(m,\varphi)}{\varphi_{+}(a)},\;\;\for\;a>0.
$$
We can obviously take $a$ such that 
$$
\frac{B\ccM(m,\varphi)}{\varphi_{+}(a)}=AR\tau_{m,\varphi}(s,t),\;\;\mbox{i.e.}\;\; a=\varphi^{-1}_{+}(\frac{B\ccM(m,\varphi)}{AR\tau_{m,\varphi}(s,t)}),
$$
thus denoting $K=ARB^{-1}$ we derive
$$
\frac{|f(s)-f(t)|}{ 2AR^3\tau_{m,\varphi}(s,t)\varphi^{-1}_{+}(\frac{\ccM(m,\varphi)}{K\tau_{\varphi,m}(s,t)})}\ls |f^{d}|^{\nu}_{\varphi_{+}}.
$$
Lemma \ref{2lem} gives the result with $C=2AR^5$.
\end{dwz}
\smallskip

\noindent
{\bf Acknowledgment} I would like to thank professor Stanislaw Kwapien and the anonymous referee for numerous remarks which
helped me to improve the paper.

\end{document}